\documentclass[11pt]{amsart}
\usepackage{geometry}                
\usepackage{csquotes}
\usepackage{graphicx}
\usepackage{amssymb}
\usepackage{epstopdf}
\title[EnKBF for NMPC]{Ensemble Kalman--Bucy filtering for nonlinear model predictive control}
\author{Sebastian Reich}
\address{Institut f\"ur Mathematik, Universit\"at Potsdam, Karl-Liebknecht-Str. 24/25, 14476 Potsdam}
\date{\today}                                           

\begin{document}

\begin{abstract}
We consider the problem of optimal control for partially observed dynamical systems. Despite its prevalence in practical applications, there are still very few algorithms available, which take  uncertainties in the current state estimates and future observations into account. In other words, most current approaches separate state estimation from the optimal control problem. In this paper, we extend the popular ensemble Kalman filter to receding horizon optimal control problems in the spirit of nonlinear model predictive control. We provide an interacting particle approximation to the forward-backward stochastic differential equations arising from Pontryagin's maximum principle with the forward stochastic differential equation provided by the time-continuous ensemble Kalman--Bucy filter equations. The receding horizon control laws are approximated as linear and are continuously updated as in nonlinear model predictive control. We illustrate the performance of the proposed methodology for an inverted pendulum example.
\end{abstract}

\maketitle

%
\section{Introduction}
%

Optimal control for partially observed dynamical systems has first been consider in \cite{ASTROM1965} and theoretical foundations have been established, for example, in \cite{Bensoussan92} and \cite{Nisio15}. It is well-known \cite{Handel07,meyn22} that, for linear systems, the involved optimal control and state estimation problems can be solved separately. Such a separation or uncertainty equivalence principle does however no longer apply to nonlinear dynamical systems. Still such a principle is widely applied in practice due to the computational intractability of combined estimation and control problems.

In this paper, we revisit this problem by combining nonlinear model predictive control (NMPC) \cite{MPC} with ensemble Kalman filtering (EnKF) \cite{CRS22}; thus providing an alternative approach to currently available output feedback NMPC techniques \cite{Allgower,MPC}.

To simplify the exposition we assume a data assimilation setting where partial and noisy observations are taken continuously in time \cite{law2015data,reich2015probabilistic} and processed by the ensemble Kalman--Bucy filter (EnKBF) \cite{bergemann2012ensemble}.
The novelty of our approach arises from the fact that we use the EnKBF for (i) estimating model states from actual observations as well as (ii) predicting uncertainties in the optimal control problem arising from the statistics of future observations and uncertainties in the estimated model states. For the latter task we formulate an appropriate set of forward-backward stochastic differential equations (FBSDEs) based on Pontryagin's maximum principle \cite{Carmona}. In the spirit of the EnKBF, we also utilize a linear approximation {\it ansatz} when solving the associated backward stochastic differential equation (BSDE) giving rise to linear control laws in the estimated model states. While linear control laws become insufficient for highly nonlinear infinite horizon optimal control problems, the receding horizon approach of NMPC allows for a continuous adaptation of  linear approximations; very much in line with modern data assimilation techniques such as the EnKF. As for NMPC, the success of the proposed methodology depends on suitable choices of the optimal control cost function as well as sufficiently informative observations in order to apply EnKF-based data assimilation. 

In this paper, we assume a deterministic model dynamics. Model errors can be taken into account via stochastic perturbations in the EnKBF or via robust NPMC schemes \cite{Allgower,MPC}. 

The remainder of the paper is structured as follows. A mathematical problem formulation is provided in Section \ref{sec2}. We develop the proposed methodology in Section \ref{sec3} in the context of linear partially observed control problems first. In particular, we derive an interacting particle formulation in terms of FBSDEs which exactly reproduces the known state estimation and optimal control solution \cite{Handel07,meyn22}. The proposed methodology is then extended to nonlinear problems in Section \ref{sec4} very much in the spirit of EnKF-based methodologies \cite{CRS22} with the addition that we also need to approximate solutions to a BSDE. We introduce the acronym EnKBF-NMPC for the newly proposed method. Section \ref{sec5} provides a numerical illustration in the context of an inverted pendulum. The paper closes with some final remarks in Section \ref{sec6}.

%
\section{Mathematical problem formulation} \label{sec2}
%

We consider controlled ordinary differential equations (ODEs) of the form
\begin{equation}
\dot{X}_t = f(X_t)+ GU_t
\end{equation}
subject to uncertain initial conditions $X_0 \in \mathbb{R}^{d_x}$, which are distributed according to a known probability measure $\nu_0({\rm d}x)$. The control applied at time $t\ge 0$ is denoted by $U_t \in \mathbb{R}^{d_u}$. The drift function $f:\mathbb{R}^{d_x}\to \mathbb{R}^{d_x}$ and the matrix $G\in \mathbb{R}^{d_x\times d_u}$ are assumed to be given.

We also assume that the process is partially observed continuously in time with forward model
\begin{equation} \label{eq:obs model}
{\rm d}\mathcal{Y}_t = h(X_t){\rm d}t + R^{1/2}{\rm d}W_t,
\end{equation}
where $h:\mathbb{R}^{d_x}\to \mathbb{R}^{d_y}$ denotes the forward model, $R \in \mathbb{R}^{d_y\times d_y}$ is the error covariance matrix, and $W_t$ is standard $d_y$-dimensional Brownian motion (BM).

Actual observations are denotes by $\mathcal{Y}_t^\dagger$, $t\ge 0$, and its associated filtering distribution at time $\tau >0$ by 
\begin{equation}
\nu_\tau ({\rm d}x) := \nu_\tau({\rm d}x|\mathcal{Y}_{0:\tau}^\dagger).
\end{equation}
Here we have introduced the notation $\mathcal{Y}^\dagger_{0:\tau}$ to indicate the dependence on observations $\mathcal{Y}_t^\dagger$ from the time interval $t\in (0,\tau]$.

At each time instance $\tau_n = n \Delta \tau$, $n\ge 0$ and for available observations $\mathcal{Y}_{0:\tau_n}^\dagger$, we solve a finite-horizon stochastic optimal control problem with cost function
\begin{equation} \label{eq:gcost}
    \mathcal{L}(U_{\tau_n:\tau_n + T}) = \mathbb{E}\left[ 
    \int_{\tau_n}^{\tau_n + T} \left(\frac{1}{2} \|U_t\|^2 + c(X_t)\right) 
    {\rm d}t + \psi_T(X_{\tau_n +T})\right],
\end{equation}
$T > \tau_{n+1}$, with expectation taken over initial conditions 
\begin{equation}
X_{\tau_n} \sim \nu_{\tau_n}
\end{equation}
and all possible future observation paths $\mathcal{Y}_{\tau_n:\tau_n+T}$ satisfying the model (\ref{eq:obs model}).

The desired optimal control, denoted by $U_t^\ast$, is a function of distributions $\nu ({\rm d}x)$, which we denote by $u_t^\ast(\nu)$. In other words, one applies $U_t^\ast = u_t^\ast(\nu_t)$ over the time interval $t \in (\tau_n,\tau_{n+1}]$ along the filtering distribution $\nu_t$ for given data $\mathcal{Y}_{0:t}^\dagger$. Hence, we obtain the controlled ODE
\begin{equation}
    \dot{X}_t = f(X_t)+ Gu_t^\ast (\nu_t)
\end{equation}
for $t \in (\tau_n,\tau_{n+1}]$ and initial $X_{\tau_n} \sim \nu_{\tau_n}$.

Receding horizon optimal control formulations are at the heart of NMPC \cite{MPC}. NMPC typically replaces the uncertain $X_{\tau_n}$ by its expectation value under the filtering distribution $\nu_{\tau_n}$; i.e.~by
\begin{equation}
    \bar X_{\tau_n} := \mathbb{E}[X_{\tau_n}].
\end{equation}
This approximation reduces the stochastic optimal control problem to a standard deterministic optimal control problem with cost
\begin{equation}
\mathcal{L}(u_{\tau_n:\tau_n + T}) = 
    \int_{\tau_n}^{\tau_n + T} \left( \frac{1}{2} \|U_t\|^2 + c(\bar x_t)\right) {\rm d}t + \psi_T(\bar x_{\tau_n +T})
\end{equation}
subject to the ODE
\begin{equation}
\dot{\bar x}_t = f(\bar x_t) + GU_t
\end{equation}
with initial condition $\bar x_{\tau_n} = \bar X_{\tau_n}$. We note that the optimal control $U_t^\ast$ becomes open loop.

In this paper, we develop an alternative approach which takes the uncertainties in the initial conditions and future observations into account and results in linear (closed loop) control laws of the form
\begin{equation} \label{eq:optimal control law}
    u_t^\ast(\nu_t) = -G^{\rm T} \left( \Lambda_t^\ast \bar X_t + \lambda_t^\ast \right)
\end{equation}
for suitable matrices $\Lambda_t^\ast \in \mathbb{R}^{d_x\times d_x}$ and vectors $\lambda_t^\ast \in \mathbb{R}^{d_x}$, $t \in (\tau_{t_n},\tau_{t_n}+T]$. Our approach is based on the continuous-time EnKBF approximation \cite{CRS22} to the filtering distributions $\nu_t$ and we call our approach therefore EnKBF-NMPC. 

Let us embed the discussion so far into a digital twin setting \cite{KPW21}. The computed optimal control
$U_t^\dagger = u_t^\ast(\nu_t)$; i.e.~(\ref{eq:optimal control law}), would be applied to the physical twin, which could be of the form
\begin{equation} \label{eq:PT}
    {\rm d}X_t^\dagger = f(X_t^\dagger){\rm d}t + GU_t^\dagger {\rm d}t + \sigma {\rm d}B_t^\dagger .
\end{equation}
Here $B_t^\dagger$ stands for random perturbations, which mimic the difference between the physical and digital twin and which could take the form of standard BM for simplicity of exposition, and $\sigma > 0$. The physical twin (\ref{eq:PT}) produces a reference trajectory, which we denote by
$X_t^\dagger$ and which is partially observed according to the observation model (\ref{eq:obs model}); i.e., the actual observations $\mathcal{Y}_{0:t}^\dagger$ satisfy 
\begin{equation}
{\rm d}\mathcal{Y}_t^\dagger = 
h(X_t^\dagger){\rm d}t + {\rm d} W_t^\dagger
\end{equation}
with $W_t^\dagger$ a particular realization of BM. Given the data, the digital twin evolves according to the EnKBF equation
\begin{equation}
    {\rm d}X_t = f(X_t){\rm d}t - GG^{\rm T}\left(\Lambda_t^\ast \bar X_t + \lambda_t^\ast\right)
    {\rm d}t - \frac{1}{2} C_t^{xh}R^{-1}(h(X_t)-\bar h_t){\rm d}t +
    C_t^{xh}R^{-1} {\rm d}\mathcal{Y}_t^\dagger,
\end{equation}
which depends on the law $\nu_t$ of $X_t$ via its mean $\bar X_t$, the expectation value $\bar h_t$ of $h(X_t)$, and the covariance matrix $C^{xh}_t$ between $X_t$ and $h(X_t)$.

In order to implement the proposed approach, it remains to formulate an algorithm for computing the required $\Lambda_t^\ast$ and $\lambda_t^\ast$. Our approach is motivated in the following section by first considering linear control problems.

%
\section{Optimal control for partially observed linear ODEs} \label{sec3}
%

We start from the linear control problem
\begin{equation}
    \dot{X}_t = AX_t + b + GU_t
\end{equation}
on state space $\mathbb{R}^{d_x}$ subject to continuous-time observations
\begin{equation} \label{eq:obs}
{\rm d}\mathcal{Y}_t = HX_t {\rm d}t + R^{1/2} {
\rm d}W_t.
\end{equation}
An actual observation is denoted by $\mathcal{Y}_t^\dagger \in \mathbb{R}^{d_y}$ and the initial $X_0$ is Gaussian distributed with mean $m_0$ and covariance $C_0$. The EnKBF equation is given by
\begin{equation} \label{eq:LEnKBF}
{\rm d}X_t^{(i)} =
\left(AX_t^{(i)}+b + GU_t\right){\rm d}t
- \frac{1}{2} \mathcal{C}_t H^{\rm T} R^{-1} H(X_t^{(i)}+
\bar X_t){\rm d}t + \mathcal{C}_t H^{\rm T}R^{-1/2} {\rm d}\mathcal{Y}_t^\dagger,
\end{equation}
$i=1,\ldots,M$ \cite{CRS22}. Here $M>d_x$ is the ensemble size, 
\begin{equation}
\bar X_t^M = \frac{1}{M} \sum_{j=1}^M X_t^{(j)}
\end{equation}
denotes the empirical mean of the ensemble  $\{X_t^{(i)}\}$ and 
\begin{equation}
\mathcal{C}_t = \frac{1}{M} \sum_{j=1}^M (X_t^{(j)}-\bar X_t^M)(X_t^{(j)}-\bar X_t^M)^{\rm T}
\end{equation}
its empirical covariance matrix.

The initial ensemble is chosen such that  $\bar X_0^M=m_0$ and $\mathcal{C}_0 = C_0$. Given these two conditions, the EnKBF formulation (\ref{eq:LEnKBF}) is equivalent to the standard Kalman--Bucy filter equations in terms of the time evolution of their mean and covariance matrix; i.e., $\bar X_t^M = m_t$ and $\mathcal{C}_t = C_t$ for all $t\ge 0$ regardless of the ensemble size $M>d_x$. Note, however, that the EnKBF is often used with $M < d_x$ in combination with localization and that the EnKBF easily generalizes to non-Gaussian and nonlinear evolution settings.

We set $\tau_n=0$ in (\ref{eq:gcost}) for simplicity of exposition and consider the time window $[0,T]$ with the optimal control $U_t^\ast$ now chosen to minimize the cost functional
\begin{equation} \label{eq:cost}
\mathcal{V}(U_{0:T}) = \mathbb{E}\left[ 
\int_0^T
\left( \frac{1}{2} \|U_t\|^2 + 
\frac{1}{M} \sum_{j=1}^M c(X_t^{(j)})
 \right) {\rm d}t  +
\frac{1}{2M}\sum_{j=1}^M \psi(X_T^{(j)})\right],
\end{equation}
where
\begin{equation}
c(x) = \frac{1}{2}(x-c)^{\rm T}V(x-c)
\end{equation}
is the running and
\begin{equation}
\psi(x) = \frac{1}{2}(x-c_T)^{\rm T}V_T(x-c_T)
\end{equation}
the terminal cost. The symmetric positive-definite matrices $V \in \mathbb{R}^{d_x \times d_x}$ and $V_T\in \mathbb{R}^{d_x\times d_x}$ as well as the vectors $c \in \mathbb{R}^{d_x}$ and $c_T \in \mathbb{R}^{d_x}$ are all assumed to be given. Expectation in (\ref{eq:cost}) is taken over all possible future observations $\mathcal{Y}_t$, $t>0$.

In order to apply the stochastic Pontryagin maximum principle \cite{Carmona}, we first rewrite the EnKBF (\ref{eq:LEnKBF}) in the form
\begin{equation} \label{eq:LEnKBFn}
{\rm d}X_t^{(i)} =
\left(AX_t^{(i)}+b + GU_t\right){\rm d}t
- \frac{1}{2} \mathcal{C}_t H^{\rm T} R^{-1} H(X_t^{(i)}-
\bar X_t^M){\rm d}t + \mathcal{C}_t H^{\rm T}R^{-1/2} {\rm d}\mathcal{W}_t
\end{equation}
$t>0$, where $\mathcal{W}_t$ is standard $d_y$-dimensional BM. This reformulation is justified because the innovation process $\mathcal{I}_t$ defined by ${\rm d}\mathcal{I}_t = {\rm d}\mathcal{Y}_t- H\bar X_t{\rm d}t$ behaves like BM \cite{SR-crisan}. It should be noted that the noise $\mathcal{W}_t$ only affects the time evolution of the mean $\bar X_t^M$. 

Let us assume for simplicity that observations $\mathcal{Y}_t$ in (\ref{eq:obs}) are scalar-valued, i.e., $d_y=1$. Following \cite{Carmona}, the Hamiltonian appearing in Pontryagin's maximum principle is given by
\begin{subequations}
    \begin{align}
    & \mathcal{H}(\{x_j\},\{y_j\},\{z_j\},u) \,= \frac{M}{2} \|u\|^2 +
    \sum_{j=1}^M c(x_j)\,+\\
& \qquad \quad \sum_{j=1}^M \left\{
\left(Ax_j+b + Gu -\frac{1}{2}
\mathcal{C}H^{\rm T} R^{-1} H(x_j-\bar x)\right)^{\rm T}y_j +  R^{-1/2} H\mathcal{C}z_j \right\}
\end{align}
\end{subequations}
and the associated FBSDEs become 
\begin{subequations} \label{eq:LFBSDE}
\begin{align}
{\rm d} X_t^{(i)} &= \nabla_{y_i} H(\{X_t^{(j)}\},\{Y_t^{(j)}\},\{Z_t^{(j)}\},U_t^\ast){\rm d}t +
\mathcal{C}_t H^{\rm T} R^{-1/2} {\rm d}\mathcal{W}_t,\\
{\rm d}Y^{(i)}_t &=
-\nabla_{x_i} H(\{X_t^{(j)}\},\{Y_t^{(j)}\},\{Z_t^{(j)}\},U_t^\ast){\rm d}t + Z_t^{(i)} {\rm d}\mathcal{W}_t
\end{align}
\end{subequations}
for $i=1,\ldots,M$, subject to the terminal condition
\begin{equation}
Y_T^{(i)} = \nabla_x \psi(X_T^{(i)})
\end{equation}
and the initial ensemble $\{X_0^{(j)}\}$. The desired optimal control $U_t^\ast$ is provided by
\begin{equation} \label{eq:OC law}
    U_t^\ast = \arg \min_u H(\{X_t^{(j)}\},\{Y_t^{(j)}\},\{Z_t^{(j)}\},u) 
    = -\frac{1}{M} \sum_{j=1}^M G^{\rm T} Y_t^{(j)}.
\end{equation}

The computation of $\nabla_{x_i}\mathcal{H}$ is rather complex due to the appearance of the empirical covariance matrix $\mathcal{C}$ and the empirical mean 
\begin{equation}
\bar x = \frac{1}{M}\sum_{j=1}^M x_j
\end{equation}
in the Hamiltonian $\mathcal{H}$:
\begin{subequations} \label{eq:gradH}
    \begin{align}
    \nabla_{x_i} \mathcal{H} &= A^{\rm T} y_i - \frac{1}{2} H^{\rm T} R^{-1} H \mathcal{C}
    (y_i - \bar y) \,+\\
    & \quad - \,\frac{1}{2} H^{\rm T} R^{-1} H \mathcal{C}^{xy}(x_i-\bar x) - \frac{1}{2}
    \mathcal{C}^{yx}H^{\rm T}R^{-1} H(x_i-\bar x)\,+\\
    & \quad +\, \frac{1}{M}\sum_{j=1}^M H^{\rm T} R^{-1/2} z_j^{\rm T}(x_i-\bar x) +
    \frac{1}{M} \sum_{j=1}^M z_j R^{-1/2} H (x_i-\bar x).
    \end{align}
\end{subequations}
Here we have also introduced the empirical covariance matrix
\begin{equation}
\mathcal{C}^{yx} = \frac{1}{M} \sum_{j=1}^M (y_j-\bar y)(x_j -\bar x)^{\rm T}
\end{equation}
and the mean
\begin{equation}
    \bar y = \frac{1}{M}\sum_{j=1}^M y_j.
\end{equation}
We next introduce the ensemble deviations
\begin{equation}
\delta X_t^{(i)} := X_t^{(i)}-\bar X_t^M
\end{equation}
for $i=1,\ldots,M$ and make the {\it ansatz}
\begin{equation} \label{eq:ansatz}
Y_t^{(i)} = \Lambda_t \bar X_t^M + \Omega_t \delta X_t^{(i)} + \lambda_t
\end{equation}
for suitable matrices $\Lambda_t \in \mathbb{R}^{d_x\times d_x}, 
\Omega_t \in \mathbb{R}^{d_x\times d_x}$ and vectors $\lambda_t \in \mathbb{R}^{d_x}$. It then follows from \cite{Carmona} eq.~(4.25) that
\begin{equation}\label{eq:Z}
    Z_t^{(j)} = \Lambda_t \mathcal{C}_t H^{\rm T}R^{-1/2} 
\end{equation}
for all $j=1,\ldots,M$. Also note that 
\begin{equation}
\mathcal{C}_t^{yx} = \Omega_t \mathcal{C}_t.
\end{equation}

The optimal control is defined in terms of the ensemble mean
\begin{equation}
    \bar Y_t = \frac{1}{M}\sum_{j=1}^M Y_t^{(j)} = \Lambda_t \bar X_t^M + \lambda_t
\end{equation}
and, therefore, we only consider the evolution equation in $\bar Y_t$ from now on, which is given by the BSDE
\begin{equation} \label{eq:LBSDE}
    -{\rm d}\bar Y_t =
    A^{\rm T} \bar Y_t{\rm d}t  + V (\bar X_t^M - c){\rm d}t - \bar Z_t {\rm d}\mathcal{W}_t.
\end{equation}
Furthermore, the forward stochastic differential equation (FSDE) (\ref{eq:LFBSDE}a) becomes
\begin{equation} \label{eq:LFSDE}
{\rm d}X_t^{(i)} =
\left(AX_t^{(i)}+b + G U_t^\ast\right){\rm d}t
- \frac{1}{2} \mathcal{C}_t H^{\rm T}R^{-1} H(X_t^{(i)}-
\bar X_t^M){\rm d}t + \mathcal{C}_t H^{\rm T}R^{-1/2} {\rm d}\mathcal{W}_t,
\end{equation}
with optimal control law 
\begin{equation}
U_t^\ast = u_t^\ast(\nu_t):= - G^{\rm T}\bar Y_t = -G^{\rm T}
\left( \Lambda_t \bar X_t^M + \lambda_t\right).
\end{equation}

Hence (\ref{eq:LFSDE}) and (\ref{eq:LBSDE}) together with (\ref{eq:ansatz}) lead to
the following evolution equations for $\Lambda_t$ and $\lambda_t$:
\begin{subequations} \label{eq:Lambda}
    \begin{align}
    -\dot{\Lambda}_t &= \Lambda_t A + A^{\rm T}\Lambda_t - \Lambda_t G G^{\rm T} \Lambda_t 
    + V ,\\
    -\dot{\lambda}_t &= A^{\rm T} \lambda_t - V c + \Lambda_t (b - GG^{\rm T} \lambda_t)
    \end{align}
\end{subequations}
with terminal conditions $\Lambda_T = V_T$ and $\lambda_T = -V_T c_T$, respectively. Here we have taken the ensemble average which eliminates all terms involving
$\Omega_t$ and $\delta X_t^{(j)}$, $j=1,\ldots,M$.

The formulation (\ref{eq:Lambda}) coincides with the equations from optimal control for linear systems \cite{Handel07,meyn22}. We have derived them here from a  perspective based on the McKean (interacting particle) reformulation (\ref{eq:LEnKBF}) of the standard Kalman--Bucy filter equations \cite{CRS22} and an application of the stochastic Pontryagin maximum principle \cite{Carmona}.

%
\section{The EnKBF-NMPC algorithm} \label{sec4}
%

The attractiveness of formulations (\ref{eq:LEnKBFn}) and (\ref{eq:LBSDE}) comes from the fact that they generalize beyond the linear Gaussian setting. We now also consider vector-valued observations; i.e.~$d_y\ge 1$. For example, a generalized FBSDE formulation is provided by
\begin{subequations} \label{eq:FBSDE}
    \begin{align}
    {\rm d}X_t^{(i)} &=
(f(X_t^{(i)})+GU_t^\ast){\rm d}t
- \frac{1}{2} \mathcal{C}_t^{xh} R^{-1} (h(X_t^{(i)})- \bar h_t^M
){\rm d}t - \mathcal{C}^{xh}_t R^{-1/2} {\rm d}\mathcal{W}_t,\\
-{\rm d}\bar Y_t &=
    Df(\bar X_t^M)^{\rm T} \bar Y_t{\rm d}t  +
     \nabla_x c(\bar X_t^M){\rm d}t +  \bar Z_t {\rm d} \mathcal{W}_t,
    \end{align}
\end{subequations}
where
\begin{equation}
\bar h_t^M = \frac{1}{M}\sum_{j=1}^M h(X_t^{(j)})
\end{equation}
and $\mathcal{C}_t^{xh}$ denotes empirical covariance matrices between $X_t \in \mathbb{R}^{d_x}$ and $h(X_t) \in \mathbb{R}^{d_y}$, respectively. The desired optimal control is still provided by
\begin{equation}
    U_t^\ast = u_t^\ast (\nu_t) := - G^{\rm T}\bar Y_t.
\end{equation}

Formulation (\ref{eq:FBSDE}) is now used as as the basis of extending NMPC \cite{MPC} to partially observed diffusion processes and EnKBF-based data assimilation. In particular, we have to simulate (\ref{eq:FBSDE}) over multiple realizations of the noise $\mathcal{W}_t$, which we denote by $\mathcal{W}_{t,k}$, $k=1,\ldots,K$. We also apply statistical linearization \cite{CRS22} in order to avoid computation of the Jacobian $Df(x)$. These approximations yield
\begin{subequations} \label{eq:SFBSDE}
    \begin{align}
    {\rm d}X_{t,k}^{(i)} &= \left(f(X_{t,k}^{(i)})+GU_{t,k}^\ast
- \frac{1}{2} \mathcal{C}_{t,k}^{xh} R^{-1}(h(X_{t,k}^{(i)})- \bar h^M_{t,k}
)\right){\rm d}t + \mathcal{C}^{xh}_{t,k}R^{-1/2} {\rm d}\mathcal{W}_{t,k},\\
-{\rm d}\bar Y_{t,k} &=
    \mathcal{C}_{t,k}^{-1}\mathcal{C}_{t,k}^{xf} \bar Y_{t,k}{\rm d}t
 +  \nabla_x c(\bar X_{t,k}^M){\rm d}t - \bar Z_{t,k} {\rm d} \mathcal{W}_{t,k},
    \end{align}
\end{subequations}
$k=1,\ldots,K$, together with optimal control $U_{t,k}^\ast = -G^{\rm T} \bar Y_{t,k}$. The empirical covariance matrices are defined by
\begin{equation}
\mathcal{C}^{xh}_{t,k} = \frac{1}{M}\sum_{j=1}^M (X_{t,k}^{(j)}-\bar X_{t,k}^M)(
h(X_{t,k}^{(j)})-\bar h^M_{t,k})^{\rm T}
\end{equation}
and related expressions for $\mathcal{C}_{t,k} = \mathcal{C}_{t,k}^{xx}$ and
$\mathcal{C}_{t,k}^{xf}$. Furthermore,
\begin{equation}
\bar X_{t,k}^M = \frac{1}{M} \sum_{j=1}^M X_{t,k}^{(j)},
 \qquad
\bar h_{t,k}^M = \frac{1}{M}\sum_{j=1}^M h(X_{t,k}^{(j)}).
\end{equation}

The BSDE (\ref{eq:SFBSDE}b) is now approximated using the linear regression {\it ansatz}
\begin{equation}
\bar Y_{t,k} \approx \Lambda_t (\bar X_{t,k}^M-\bar X_t^M) + \mu_t,
\end{equation}
where $\bar \Lambda_t$, $\mu_t$, and $\bar Z_{t,k}$ are chosen such that the resulting approximation error is minimized \cite{Carmona} and
\begin{equation}
 \bar X_t^M = \frac{1}{K} \sum_{k=1}^K \bar X_{t,k}^M.
\end{equation}
We note that the choice
\begin{equation} \label{eq:open loop MPC}
    U^\ast_{t,k} = -G^{\rm T} \mu_t
\end{equation}
corresponds to an open loop control in NMPC.

More precisely, given a time-discretization at time-levels 
$t_{n+1} = t_n+\Delta t$, $n=0,\ldots,N-1$
for given step-size $\Delta t = T/N>0$ and $t_0 = 0$, we define the least squares cost function
\begin{equation} \label{eq:cost 1}
\mathcal{L}_n(\Lambda_{t_n},\mu_{t_n}) = \frac{1}{K} \sum_{k=1}^K
\left\|\tilde Y_{t_n,k}
- (I + \Delta t \mathcal{C}_{t_{n+1},k}^{-1}
\mathcal{C}_{t_{n+1},k}^{xf})\bar Y_{t_{n+1},k} -
\nabla_x c(\bar X_{t_{n+1},k}^M)\Delta t \right\|^2
\end{equation}
subject to 
\begin{equation} \label{eq:linear regression 1}
\tilde Y_{t_n,k} = \Lambda_{t_n}(\bar X_{t_n,k}^M-\bar X_{t_n}^M)+ \mu_{t_n}.
\end{equation}

We further simplify this formulation by explicitly removing the noise in
$\bar X_{t_n,k}^M$ using the back integration step
\begin{equation} \label{eq:BODE}
    \tilde X_{t_n,k} = \bar X_{t_{n+1},k}^M + \Delta t \left( f(\bar X_{t_{n+1},k}^M) -
    G G^{\rm T}\bar Y_{t_{n+1},k} \right)
\end{equation}
for given $\bar X_{t_{n+1},k}^M$ and $\bar Y_{t_{n+1},k}$. One then replaces (\ref{eq:linear regression 1}) by
\begin{equation}
\tilde Y_{t_n,k} = \Lambda_{t_n}(\tilde X_{t_n,k}-\tilde X_{t_n})+ \mu_{t_n},\qquad \tilde X_{t_n} = \frac{1}{K}\sum_{k=1}^K \tilde X_{t_n,k},
\end{equation}
and (\ref{eq:cost 1}) becomes
\begin{equation} \label{eq:cost 2}
\mathcal{L}_n(\Lambda_{t_n},\mu_{t_n}) = \frac{1}{K} \sum_{k=1}^K
\left\|\Lambda_{t_n}(\tilde X_{t_n,k}-\tilde X_{t_n})+ \mu_{t_n} - \gamma_{t_n,k} \right\|^2
\end{equation}
with
\begin{equation}
\gamma_{t_n,k} = (I + \Delta t \mathcal{C}_{t_{n+1},k}^{-1}
\mathcal{C}_{t_{n+1},k}^{xf})\bar Y_{t_{n+1},k} +
\nabla_x c(\bar X_{t_{n+1},k}^M)\Delta t.
\end{equation}

Denote the minimizer of (\ref{eq:cost 2}) by $(\Lambda^\ast_{t_n},\mu_{t_n}^\ast)$. See the Appendix for explicit solution formulas. Then we set
\begin{equation} 
    \bar Y_{t_n,k} := \Lambda_{t_n}^\ast (\bar X_{t_n,k}^M-\tilde X_{t_n}) +  \mu^\ast_{t_n}=
    \Lambda_{t_n}^\ast \bar X_{t_n,k}^M + \lambda^\ast_{t_n}
\end{equation}
with
\begin{equation}
    \lambda^\ast_{t_n} = \mu^\ast_{t_n} - \Lambda_{t_n}^\ast\tilde X_{t_n}.
\end{equation}
The iteration is started at time $t_N = T$ with terminal condition provided by
\begin{equation}
\bar Y_{T,k} = \nabla_x \psi(\bar X_{T,k}^M).
\end{equation}

Formulation (\ref{eq:BODE}) is justified in our setting since we assume linear McKean control laws and hence second-order terms, which would arise from the stochastic contributions, do not appear in the control. 
We also note that $\Lambda_t$ is a symmetric matrix in the linear setting which is no longer the case when used as an approximation in the nonlinear setting. A symmetric approximation could still be obtained by using
$(\Lambda_{t_n}^\ast + (\Lambda_{t_n}^\ast)^{\rm T})/2$ instead of $\Lambda_{t_n}^\ast$ in (\ref{eq:control approximation}) below.

The sequence $\{\Lambda_{t_n}^\ast,\lambda_{t_n}^\ast\}_{n=0}^N$ is stored and linearly interpolated to provide the linear control 
\begin{equation} \label{eq:control approximation}
U_t^\ast = u_t^\ast(\nu_t) := -G^{\rm T}\left( \Lambda_t^\ast \bar X_t^M + \lambda_t^\ast \right)
\end{equation}
for the physical twin (\ref{eq:PT}) as well as for the forward prediction step of its digital twin
\begin{equation} \label{eq:controlled EnKBF}
    {\rm d}X_{t}^{(i)} = \left(f(X_{t}^{(i)})+GU_t^\ast
- \frac{1}{2} \mathcal{C}_{t}^{xh} R^{-1}(h(X_{t}^{(i)})+ \bar h_{t}^M
)\right){\rm d}t + \mathcal{C}^{xh}_{t}R^{-1/2} {\rm d}\mathcal{Y}_{t}^\dagger.
\end{equation}

Let us briefly comment on the relation of the proposed EnKBF-NMPC scheme to available NMPC for partially observed processes. NMPC for partially observed processes is called output NMPC in Section 4 of \cite{MPC} and output feedback NMPC in \cite{Allgower}. Most existing approaches rely on what is called the nominal solution in \cite{MPC}; which in our context amounts to setting $W_{t}\equiv 0$ and $X_0 = \bar X_0$ and thus reduced the associated optimal control problem to a deterministic finite-horizon control problem or, alternatively, to an open loop control of the form (\ref{eq:open loop MPC}). Instead, EnKBF-NMPC takes uncertainties into account which are induced by uncertain initial conditions as well as uncertain future observations.

%
\section{Numerical example} \label{sec5}
%

\begin{figure}
\includegraphics[width=0.6\textwidth]{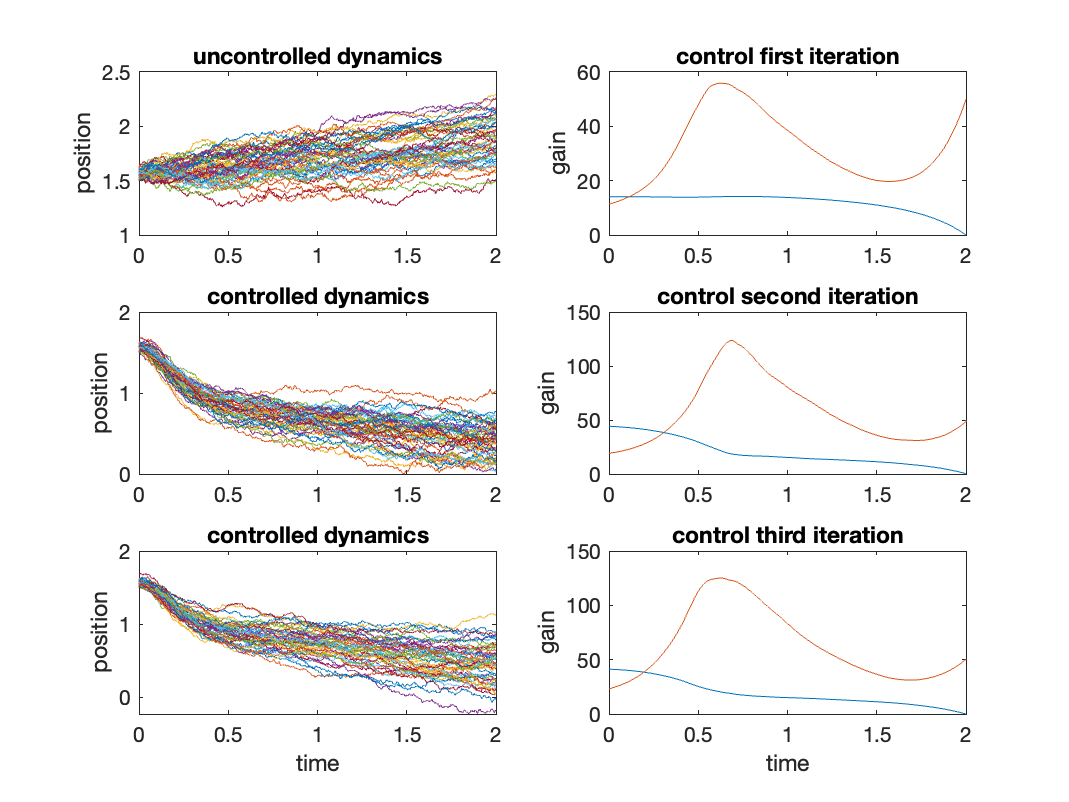}
\caption{Left panels: We display the frequentist forecast uncertainties obtained from simulating (\ref{eq:SFBSDE}a) using the optimal control obtained from the previous iteration of (\ref{eq:SFBSDE}) with the initial control being $u\equiv 0$. We perform three iterations. Right panels: We display the corresponding entries of $G^{\rm T}\Lambda_t^\ast$, which have been obtained by minimizing (\ref{eq:cost 2}). We find that the control terms have essentially converged after two iteration.} 
\label{fig1}
\end{figure}

A simple nonlinear example is provided by the controlled pendulum
\begin{equation}
\ddot{\phi}_t = \sin \phi_t -\gamma \dot{\phi}_t + u_t
\end{equation}
with $\gamma = 5$. We note that the equilibrium $(0,0)$ is unstable while $(\pi,0)$ is stable. The control $u_t$ should stabilize the unstable equilibrium given uncertainty in the initial conditions and noisy observations of the angle $\phi_t$; i.e.,
\begin{equation}
    {\rm d}\mathcal{Y}_t = \phi_t{\rm d}t + R^{1/2}{\rm d}W_t.
\end{equation}
The initial positions $\phi_0$ are Gaussian distributed with mean $\pi/2$ and the initial velocities $\dot{\phi}_0$ are Gaussian distributed with mean zero. The variance is $\sigma_0 = 0.1$ for both variables. We set $R=1$ unless stated otherwise. Observations $\mathcal{Y}_t^\dagger$, $t\ge 0$, are generated under an ideal twin experimental setting; i.e., we set $\sigma=0$ in the corresponding physical twin (\ref{eq:PT}).

We first implement a fixed-horizon version of the EnKBF-NMPC subject to final cost
\begin{equation}
\psi_T(\phi,\dot{\phi}) = \frac{\lambda_T}{2} (\phi^2+\dot{\phi}^2)
\end{equation}
and running cost
\begin{equation}
    c(\phi,\dot{\phi}) = \frac{\lambda}{2}(\phi^2 + \dot{\phi}^2)
\end{equation}
with $\lambda_T = \lambda = 50$. The implemented EnKBF-NMPC scheme starts from $u_t \equiv 0$ in the FSDE and then iterate forward and backward three times over the fixed time interval $[0,2]$; i.e.~$T=2$. The results can be found in Figure \ref{fig1}. The ensemble sizes are $M=K=50$ and the step-size is set to $\Delta t = 10^{-3}$. The effect of the control can be clearly seen and also the (frequentist) uncertainty due to possible future observations. 

\begin{figure}
\includegraphics[width=0.6\textwidth]{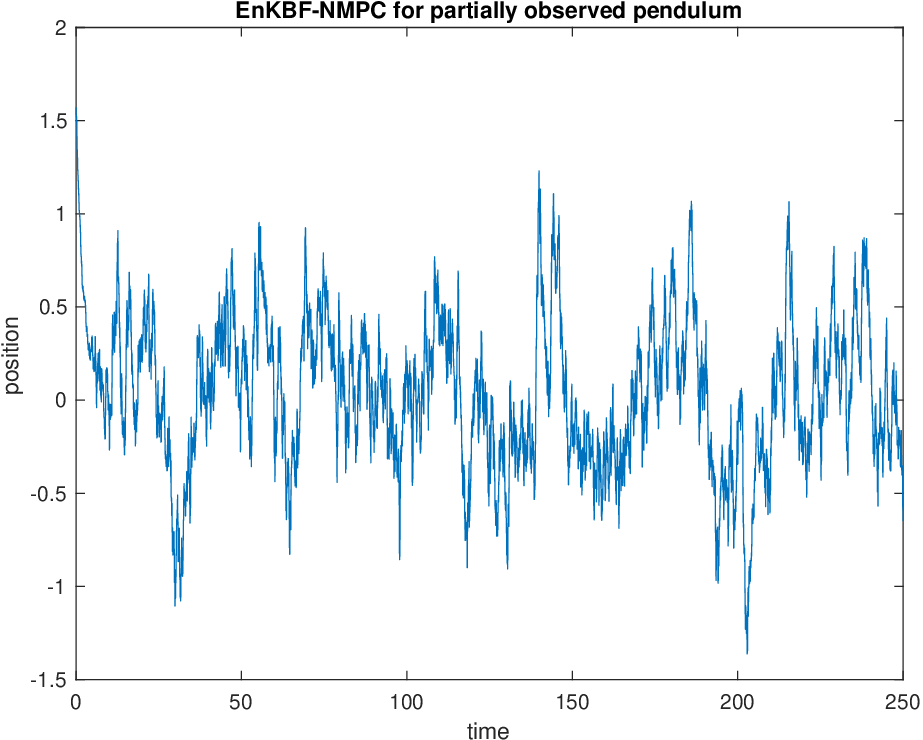}
\caption{Ensemble mean from an EnKBF filter implementation with controls computed using the proposed EnKBF-NMPC methodology. The control laws were computed over time intervals of length $T = 0.5$ with the computed control then applied over intervals of length $\Delta \tau = 0.05$. The process was repeated 5000 times. One can clearly identify the stabilizing effect of the control.} 
\label{fig2}
\end{figure}

\begin{figure}
\includegraphics[width=0.6\textwidth]{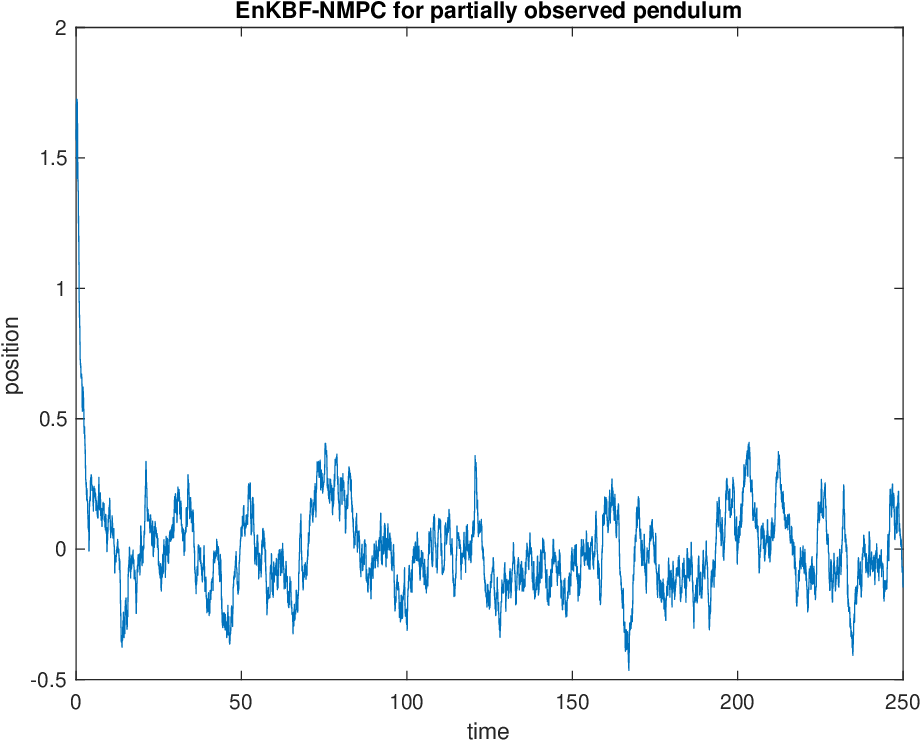}
\caption{Same setting as used for Figure \ref{fig2} except for a reduced measurement with $R = 0.1$. The reduction in variance of the computed solution relative to the unstable $q=0$ equilibrium is clearly visible.} 
\label{fig3}
\end{figure}

In a second experiment, we repeatedly solved the receding NMPC problem over time intervals of length $T = 0.5$ with final cost $\psi_{T}(x)$ and running cost $c(x)$ as stated before. The computed control is then applied over time intervals of length $\Delta \tau = 0.05$ subject to a noisy observation of the angle from its physical twin.  The resulting time evolution of the ensemble mean from the associated controlled EnKBF (\ref{eq:controlled EnKBF}) with the whole process repeated 5000 times can be found for $R=1$ in Figure \ref{fig2} and for $R=0.1$ in Figure \ref{fig3}.

%
\section{Conclusions} \label{sec6}
%

We have proposed a combination of NMPC and the EnKBF in order to approximate the effect of future observations and uncertainties in current state estimates on the control of partially observed dynamical systems. 

While this paper has taken a purely algorithmic approach, it remains to investigate the stability of the EnKBF-NMPC method under idealized modeling scenarios and its dependence on the chosen terminal cost function $\psi_T$. See, for example, \cite{de2018long} for a related accuracy and stability study for the EnKBF alone. 

Another possible extension is to replace the EnKBF by the feedback particle filter \cite{MehtaMeyn13}. Note that \cite{MehtaMeyn13} relies on the value function and its associated Hamilton--Jacobi--Bellman equation approach while this paper follows Pontryagin's maximum principle.

It would be of interest to also investigate CVaR risk measures and their sample approximations as cost function; i.e., replace $\mu_T[\psi]$ by
\begin{equation}
\tilde \psi_\beta (\mu_T) := \arg \min_s\{ s + \frac{1}{1-\beta} \mu_T \left[ (\psi_T - s)_+ \right]\} 
\end{equation}
for given $\beta \in (0,1)$. Alternatively, one could consider a terminal cost
\begin{equation}
    \tilde \psi_\beta(\mu_T) = 
    \frac{1}{\beta}\log \left(\mu_T\left[\exp(\beta \psi_T )\right] \right),
\end{equation}
for $\beta \in \mathbb{R}$ \cite{Carmona}.

We finally mention that our approach can be extended to discrete-time observations along the lines of \cite{BDRT24}; again replacing the Hamilton--Jacobi--Bellman equation by Pontryagin's maximum principle.


\medskip \medskip 

\paragraph{Acknowledgments.}

This work has been partially funded by Deutsche Forschungsgemeinschaft (DFG) - Project-ID 318763901 - SFB1294.  

\bibliographystyle{plainurl}
%
\bibliography{bib-database}
%

%
%

%
\section*{Appendix}
%

The minimizer of a functional of the form
\begin{equation}
    \mathcal{L}(\Lambda,\mu) = \sum_{k=1}^K 
    \| \Lambda (X_k-\bar X) + \mu - \gamma_k\|^2, \qquad \bar X = \frac{1}{K}\sum_{k=1}^M X_k,
\end{equation}
is provided by
\begin{equation}
    \mu^\ast = \frac{1}{K}\sum_{k=1}^K \gamma_k
\end{equation}
and
\begin{subequations}
    \begin{align}
    \Lambda^\ast &= \left(\frac{1}{K} 
    \sum_{k=1}^K (\gamma_k-\mu^\ast) (X_k-\bar X)^{\rm T}\right)
    \left( \frac{1}{K} \sum_{k=1}^K 
    (X_k-\bar X)(X_k-\bar X)^{\rm T}\right)^{-1} \\
    &= \mathcal{C}^{\gamma x}(\mathcal{C}^{xx})^{-1},
\end{align}
\end{subequations}
which implies a form of statistical linearization as also commonly used in ensemble Kalman filtering \cite{CRS22}.

\end{document}